\documentclass[12pt]{article}
\usepackage{amsfonts,latexsym,rawfonts,amsmath,amssymb,amsthm,graphicx}
\numberwithin{equation}{section}
\newtheorem{theorem}{Theorem}[section]

\newtheorem{thm}[theorem]{Theorem}
\newtheorem{pro}[theorem]{Proposition}

\newtheorem{rem}[theorem]{Remark}
\def\s{\,\,\,\,}

\def\endproof{$\hfill\Box$\\}

\title{On $n$-dimensional complete self-similar solutions to the mean curvature flow in $\mathbb{R}^{n+1}$ with nonnegative constant scalar curvature}
\author{Yong Luo, Linlin Sun, Jiabin Yin}
\date{}
\begin{document}
\maketitle
\begin{abstract}
As is well known, self-similar solutions to the mean curvature flow, including self-shrinkers, translating solitons and self-expanders, arise naturally in the singularity analysis of the mean curvature flow. Recently, Guo \cite{Guo} proved that $n$-dimensional compact self-shrinkers in $\mathbb{R}^{n+1}$ with scalar curvature bounded from above or below by some constant are isometric to the round sphere $\mathbb{S}^n(\sqrt{n})$, which implies that  $n$-dimensional compact self-shrinkers in $\mathbb{R}^{n+1}$ with constant scalar curvature are isometric to the round sphere $\mathbb{S}^n(\sqrt{n})$(see also \cite{Hui1}). Complete classifications of $n$-dimensional translating solitons in $\mathbb{R}^{n+1}$ with nonnegative constant scalar curvature and of $n$-dimensional self-expanders in $\mathbb{R}^{n+1}$ with nonnegative constant scalar curvature were given by Mart\'{i}n, Savas-Halilaj and Smoczyk\cite{MSS} and Ancari and Cheng\cite{AC}, respectively. In this paper we give complete classifications of $n$-dimensional complete self-shrinkers in $\mathbb{R}^{n+1}$ with nonnegative constant scalar curvature. We will also give alternative proofs of the classification theorems due to Mart\'{i}n, Savas-Halilaj and Smoczyk \cite{MSS} and Ancari and Cheng\cite{AC}.
\end{abstract}
\section{Introduction}
\subsection{Self-shrinkers}
An $n$-dimensional submanifold  $X:M\to \mathbb{R}^{n+p}$ in the $(n+p)$-dimensional Euclidean space $\mathbb{R}^{n+p}$ is called a self-shrinker if it satisfies
\begin{eqnarray}\label{eqmain}
H=-X^N,
\end{eqnarray}where $H$ and $X^N$ denote the mean curvature vector field and the orthogonal projection of $X$ into the normal bundle of $M^n$, respectively. It is well known that self-shrinkers play an important role in the singularity analysis of the mean curvature flow as they describe all possible Type I blow ups at a given singularity of the mean curvature flow (cf. \cite{Hui1} and \cite{Hui2}). Abresch and Langer \cite{AL} gave a complete classification of closed self-shrinkers of dimension one in $\mathbb{R}^2$, which are now called Abresch-Langer curves. In the hypersurface case, Huisken \cite{Hui1, Hui2} proved a classification theorem that the only possible smooth self-shrinkers $M^n$ in $\mathbb{R}^{n+1}$ with nonnegative mean curvature, bounded second fundamental form and polynomial volume growth are isometric to $\Gamma\times\mathbb{R}^{n-1}$ or $\mathbb{S}^k(\sqrt{k})\times\mathbb{R}^{n-k}(0\leq k\leq n)$. Here $\Gamma$ is an Abresch-Langer curve. In their pioneering work, Colding and Minicozzi \cite{CM} showed that Huisken's classification theorem still hold without the assumption of bounded second fundamental form. Since then, many interesting works on self-shrinkers are done.

 Le and Sesum \cite{LS} proved that if $M$ is an $n$-dimensional complete self-shrinker with polynomial volume growth and $S<1$, then $M$ is isometric to the hyperplane $\mathbb{R}^n,$ where $S$ denotes the squared norm of the second fundamental form. Furthermore, Cao and Li \cite{CL} have studied the general case. They proved that if $M$ is an $n$-dimensional complete self-shrinker with polynomial volume growth and $S\leq1$ in the Euclidean space $\mathbb{R}^{n+p}$, then $M$ is isometric to either the hyperplane $\mathbb{R}^{n}$, the round sphere $\mathbb{S}^n(\sqrt{n})$ or a cylinder $\mathbb{S}^m(\sqrt{m})\times \mathbb{R}^{n-m}, 1\leq m\leq n-1$. See also Ding and Xin \cite{DX}, Cheng and Wei \cite{CW}, Xu and Xu \cite{XX}, and Lei and Xu and Xu \cite{LXX} for the second gap on the squared norm of the second fundamental form for  $n$-dimensional complete self-shrinkers with polynomial volume growth in Euclidean space $\mathbb{R}^{n+1}$. These extrinsic rigidity theorems for self-shrinkers reveal the similarity of self-shrinkers in a Euclidean space with minimal submanifolds in the unit sphere.

 Bewaring in mind the isoparametric conjecture for closed minimal hypersurfaces of constant squared norm of the second fundamental form in the unit sphere (cf. \cite{GT}\cite{SWY} for a survey), it is very natural for Cheng and his collaborators to study the classification of $n$-dimensional complete self-shrinkers in $\mathbb{R}^{n+1}$ with constant squared norm of the second fundamental form. The gave a complete classification of such self-shrinkers of dimension 2 in $\mathbb{R}^3$ (cf. \cite{CO}) and of dimension 3 in $\mathbb{R}^4$ with an additional assumption of constant $f_4$ very recently (cf. \cite{CLW}). Note that for minimal submanifolds in a unit sphere, by the Gauss equation, constant squared norm of the second fundamental form implies constant scalar curvature and vice versa. But for self-shrinkers this is not true and it is natural to ask if we can give a complete classification of $n$-dimensional complete self-shrinkers in $\mathbb{R}^{n+1}$ with constant scalar curvature. In this direction, Guo proved the following result.
\begin{thm}[cf. \cite{Guo}, Corollary 3.2]\label{thmguo}
An $n$-dimensional compact self-shrinker in $\mathbb{R}^{n+1}$ with constant scalar curvature is isometric to $\mathbb{S}^n(\sqrt{n})$.
\end{thm}
Theorem \ref{thmguo} is also a corollary of Huisken's previously mentioned result that the only possible smooth self-shrinkers $M^n$ in $\mathbb{R}^{n+1}$ with nonnegative mean curvature, bounded second fundamental form and polynomial volume growth are isometric to $\Gamma\times\mathbb{R}^{n-1}$ or $\mathbb{S}^k(\sqrt{k})\times\mathbb{R}^{n-k}(0\leq k\leq n)$, since if $M$ has constant scalar curvature and compact, it has nonnegative mean curvature , bounded second fundamental form and polynomial volume growth immediately. Note that Guo\cite{Guo} actually proved a more general result, that $M$ is a round sphere if it has certain lower or upper bounded on its scalar curvature and then Theorem \ref{thmguo} is an immediate corollary.

This paper aims to classify $n$-dimensional complete self-shrinkers in $\mathbb{R}^{n+1}$ with constant scalar curvature. Guo used integral formulas, in particular Minkowski's integral formulas to prove his main theorem, which seems not applicable in the complete noncompact case.  In this paper we will use point-wise estimates to give classifications of complete self-shrinkers with constant scalar curvature. We have
\begin{thm}\label{main thm}
Let $X:M\to \mathbb{R}^{n+1}$ be an $n$-dimensional complete self-shrinker in $\mathbb{R}^{n+1}$. If the scalar curvature of $M$ is constant, positive and $M$ is of polynomial volume growth, it is isometric to one of the following:
\\\rm (1) a cylinder $\mathbb{S}^m(\sqrt{m})\times \mathbb{R}^{n-m}$, $2\leq m\leq n-1$,
\\\rm (2) the round sphere $\mathbb{S}^n(\sqrt{n})$.
\end{thm}
\begin{rem}
When $M$ is compact all of assumptions in Theorem \ref{main thm} are satisfied immediately. Therefore Theorem \ref{main thm} generalizes Corollary 3.2 in \cite{Guo}. Since our proof of Theorem \ref{main thm} uses point-wise estimates, we also provide an alternative proof of Corollary 3.2 in \cite{Guo}.
\end{rem}
The main step in the proof of Theorem \ref{main thm} is to prove that if an $n$-dimensional complete self-shrinker in $\mathbb{R}^{n+1}$ has constant positive scalar curvature, its squared norm of the second fundamental form  is smaller than or equal to 1. Then the conclusion follows from Theorem 1.1 of \cite{CL}. Cao and Li  expected that the assumption of polynomial volume growth can be removed in Theorem 1.1 of \cite{CL}. Thus we except that the condition on volume growth in Theorem \ref{main thm} can be removed. We would like to propose the following conjecture.
\\

\textbf{Conjecture:} Let $X:M\to \mathbb{R}^{n+1}$ be an $n$-dimensional complete self-shrinker in $\mathbb{R}^{n+1}$. If the scalar curvature of $M$ is constant and positive, it is isometric to  a cylinder $\mathbb{S}^m(\sqrt{m})\times \mathbb{R}^{n-m}$, $2\leq m\leq n-1$ or  the round sphere $\mathbb{S}^n(\sqrt{n})$.
\\

We will prove that this conjecture is true when $M$ has constant scalar curvature larger than $n-2$ in the appendix(see Proposition \ref{thmconj}).

For $n$-dimensional complete self-shrinkers with zero scalar curvature we have
\begin{thm}\label{main thm2}
Let $X:M\to \mathbb{R}^{n+1}$ be an $n$-dimensional complete self-shrinker in $\mathbb{R}^{n+1}$. If the scalar curvature of $M$ is zero, it is isometric to $\Gamma\times\mathbb{R}^{n-1}$, where $\Gamma$ is a complete self-shrinker in $\mathbb{R}^2$.
\end{thm}
\begin{rem}
Note that self-shrinkers in $\mathbb{R}^2$ has been complete classified by Halldorsson(see \cite{Hal}, Theorem 5.1) and the closed case was earlier classified by Abresch and Langer\cite{AL}.
\end{rem}
As a corollary of Theorem \ref{main thm2}, we have
\begin{thm}\label{main thm2'}
Let $X:M\to \mathbb{R}^{n+1}$ be an $n$-dimensional complete self-shrinker in $\mathbb{R}^{n+1}$. If the scalar curvature of $M$ is zero and $M$ is of polynomial volume growth, it is isometric to one of the following:
\\ \rm(1) $\mathbb{R}^n$,
\\\rm(2) $\mathbb{S}^1\times\mathbb{R}^{n-1}$,
\\\rm(3) $\Gamma\times\mathbb{R}^{n-1}$,
where $\Gamma$ is an Abresch-Langer curve.
\end{thm}
\subsection{Translators}
An $n$-dimensional submanifold  $X:M\to \mathbb{R}^{n+p}$ in the $(n+p)$-dimensional Euclidean space $\mathbb{R}^{n+p}$ is called a translating soliton(abbreviated by translator) if it satisfies
\begin{eqnarray}\label{eqmain}
H=V^N,
\end{eqnarray}
where $H$ and $V^N$ denote the mean curvature vector field and the orthogonal projection of $V$ into the normal bundle of $M^n$, respectively. Here $V$ is a (nonzero)constant vector in $\mathbb{R}^{n+p}$.

The translators play an important role in the study of the mean curvature flow, as they often occur as Type-II singularity of a mean curvature flow (cf. \cite{AV1, AV2}\cite{HS}\cite{Wh1, Wh2}).

It is well known that there are no closed translators. Mart\'{i}n, Savas-Halilaj and Smoczyk gave a complete classification of $n$-dimensional complete translators in $\mathbb{R}^{n+1}$ with zero scalar curvature. In this paper we will show that there exists no $n$-dimensional complete translators in $\mathbb{R}^{n+1}$ with positive constant scalar curvature and we also give a slight different proof of their theorem.
\begin{thm}[\cite{MSS}]\label{main thm2.1}
Let $X:M\to \mathbb{R}^{n+1}$ be an $n$-dimensional complete translator in $\mathbb{R}^{n+1}$. If the scalar curvature $R$ of $M$ is a nonnegative constant, then we have $R=0$ and $M$ is isometric to $\Gamma\times\mathbb{R}^{n-1}$, where $\Gamma$ is a straight line or the Grim Reaper curve in $\mathbb{R}^{2}$, i.e. $M$ is a hyperplane or a grim hyperplane.
\end{thm}

\subsection{Self-expanders}
An $n$-dimensional submanifold  $X:M\to \mathbb{R}^{n+p}$ in the $(n+p)$-dimensional Euclidean space $\mathbb{R}^{n+p}$ is called a self-expander if it satisfies
\begin{eqnarray}\label{eqmain}
H=X^N,
\end{eqnarray}
where $H$ and $X^N$ denote the mean curvature vector field and the orthogonal projection of $X$ into the normal bundle of $M^n$, respectively.

Self-expanders arise naturally when one considers solutions of graphical mean curvature flow. In the case of codimension 1 and under certain assumptions on the initial hypersurface at infinity, Ecker and Huisken\cite{EH} showed that the solutions of mean curvature flow of entire graphs in Euclidean space exist for all times $t>0$ and become asymptotically self-expanding as $t\to\infty$. As was pointed out in \cite{EH} and \cite{Sm1}, self-expanders also arise as solutions of the mean curvature flow, if the initial submanifold is a cone. Moreover, in some situations uniqueness of self-expanders is important for the construction of mean curvature flows starting from certain singular configurations \cite{BM}. The interested reader can consult the recent article \cite{Sm2}  by Smoczyk and references therein for more recent progress and a historical note on self-expanders.

It is well known that self-expanders can not be closed. Ancari and Cheng\cite{AC} gave a complete classification of $n$-dimensional self-expanders in $\mathbb{R}^{n+1}$ with nonnegative constant scalar curvature. In this paper we will give a slight different proof of it.
\begin{thm}[\cite{AC}]\label{main thm3.1}
Let $X:M\to \mathbb{R}^{n+1}$ be an $n$-dimensional complete self-expander in $\mathbb{R}^{n+1}$. If the scalar curvature $R$ of $M$ is a nonnegative constant, then $R=0$ and $M$ is isometric to $\Gamma\times \mathbb{R}^{n-1}$, where $\Gamma$ is a complete self-expander in $\mathbb{R}^2$, i.e. $M$ is a hyperplane or a self-expanding hyperplane.
\end{thm}
\begin{rem}
Note that self-expanders in $\mathbb{R}^2$ has been complete classified(cf. Ishimura\cite{Ish} and Halldorsson\cite{Hal}, Theorem 6.1).
\end{rem}
\begin{rem}
In \cite{Sm2}, any self-expander $M=\Gamma\times \mathbb{R}^{n-1}\subset\mathbb{R}^{n+1}$, where $\Gamma$ is a non-trial self-expanding curve in $\mathbb{R}^2$(i.e. not a straight line) is called a self-expanding hyperplane. Here we used his terminology.
\end{rem}

We would like to point out that hypersurfaces with constant scalar curvature in a Euclidean space is a very active and important research subject of differential geometry. In particular, Hilbert's theorem states that a complete surface of constant -1 sectional curvature can not be isometrically immersed in $\mathbb{R}^3$ and Hartman and Nirenberg\cite{HN} classified complete surfaces in $\mathbb{R}^3$ with nonnegative constant curvature, i.e. they must be planes, round spheres or cylinders. Therefore our theorems make sense only when $n\geq3$. There is also a well-know Yau's conjecture for constant scalar curvature hypersurfaces in a Euclidean space which states that compact hypersurfaces in $\mathbb{R}^{n+1}$ with constant scalar curvature must be isometric to a round sphere, which was confirmed by Ros\cite{Ros} under an additional assumption of embeddedness and by Cheng\cite{Che} under an additional assumption of locally conformally flatness. Cheng and Yau \cite{CY} considered the higher dimensional case of the theorem of Hilbert, Hartman and Nirenberg, where they proved that complete noncompact hypersurfaces in $\mathbb{R}^{n+1}$ with constant scalar curvature and nonnegative sectional curvature are isometric to the generalized cylinders.

Note that we have complete classifications of $n$-dimensional complete translators and self-expanders in $\mathbb{R}^{n+1}$ with nonnegative constant scalar curvature and of $n$-dimensional complete self-shrinkers in $\mathbb{R}^{n+1}$ with zero scalar curvature or with positive constant scalar curvature and polynomial volume growth. But our argument does not work for the case of negative constant scalar curvature and therefore we would like to propose the following problem.
\\

\textbf{Problem.} Is there any $n$-dimensional complete self-similar solution to the mean curvature flow in $\mathbb{R}^{n+1}$ with negative constant scalar curvature?
\begin{rem}
For self-expanders this problem was proposed by Ancari and Cheng in \cite{AC}.
\end{rem}

\textbf{Organization.} In section 2 we give some preliminary facts on self-similar solutions to the mean curvature flow. All of the main theorems will be proved in section 3. In the appendix we give a partial positive answer to the conjecture on page 3.

\section{Preliminaries}
In this section we give some notations and formulas. Let $X:M\to \mathbb{R}^{n+1}$ be an $n$-dimensional self-shrinker in $\mathbb{R}^{n+1}$. Let $\{e_1,...,e_n,e_{n+1}\}$ be a local orthonormal basis along $M$ with dual coframe $\{\omega_1,...,\omega_n,\omega_{n+1}\}$, such that $\{e_1,...,e_n\}$ is a local orthonormal basis of $M$ and $e_{n+1}$ is normal to $M$. Then we have $$\omega_{n+1}=0, \s \omega_{n+1i}=-\sum_{j=1}^nh_{ij}\omega_j, \s h_{ij}=h_{ji},$$ where $h_{ij}$ denotes the component of the second fundamental form of $M$. Denote by $$II=\sum_{i,j}h_{ij}\omega_i\otimes\omega_je_{n+1}$$ the second fundamental form of $M$ and $$\vec{H}=\sum_{j=1}^nh_{jj}e_{n+1}, \s H=\sum_{j=1}^nh_{jj}$$ the mean curvature vector field and the mean curvature of $M$, respectively. The Gauss equations and Codazzi equations are given by
\begin{eqnarray}
R_{ijkl}&=&h_{ik}h_{jl}-h_{il}h_{jk},
\\h_{ijk}&=&h_{ikj},\label{Codazzi}
\end{eqnarray}
where $R_{ijkl}$ is the component of curvature tensor and the covariant derivative of $h_{ij}$ is defined by
$$\sum_{k=1}^nh_{ijk}\omega_k=dh_{ij}+\sum_{k=1}^nh_{kj}\omega_{ki}+\sum_{k=1}^nh_{ik}\omega_{kj}.$$
Note that from the Gauss equations we have that
\begin{eqnarray}\label{gauss}
R=H^2-S,
\end{eqnarray}
where $R$ denotes the scalar curvature of $M$ and $S$ denotes the squared norm of the second fundamental form of $M$, i.e.
$$S=\sum_{i,j}h^2_{ij}.$$
\subsection{Self-shrinkers}
The following elliptic operator $\mathcal{L}$ was introduced by Colding and Minicozzi \cite{CM}:
\begin{eqnarray}
\mathcal{L}f=\Delta f-\langle X, \nabla f\rangle,
\end{eqnarray}
where $\Delta$ and $\nabla$ denote the Laplacian and the gradient operator on the self-shrinker, respectively and $\langle\cdot,\cdot\rangle$ denotes the standard inner product of $\mathbb{R}^{n+1}$. By a direct calculation, we have(cf. \cite{CM})
\begin{eqnarray*}\label{equ1}
\mathcal{L}h_{ij}=(1-S)h_{ij},
\\\mathcal{L}H=H(1-S),
\\ \mathcal{L}|X|^2=2(n-|X|^2),
\end{eqnarray*}
where $$|X|^2=\langle X,X\rangle.$$ Then we have
\begin{eqnarray}
\frac{1}{2}\mathcal{L}S&=&\sum_{i,j,k}h_{ijk}^2+S(1-S),\label{main equ1}
\\ \frac{1}{2}\mathcal{L}H^2&=&|\nabla H|^2+H^2(1-S),\label{main equ2}
\\ \frac{1}{2}\mathcal{L}R&=&|\nabla H|^2-\sum_{i,j,k}h_{ijk}^2+R(1-S),\label{Main equ3}
\end{eqnarray}
where $h_{ijk}$ is the component of the covariant derivative of the second fundamental form.
\subsection{Translators}
For translators in a Euclidean space, Xin \cite{Xin} introduced the following elliptic operator
\begin{eqnarray}
\mathcal{L}_{II}f=\Delta f+\langle V, \nabla f\rangle,
\end{eqnarray}
where $\Delta$ and $\nabla$ denote the Laplacian and the gradient operator on the translator, respectively
and $\langle\cdot,\cdot\rangle$ denotes the standard inner product of $\mathbb{R}^{n+1}$. Then we have the following Bochner formula for the squared norm of the second fundamental form and mean curvature respectively(cf. \cite{Xin}\cite{WXZ}).
\begin{eqnarray}
\frac{1}{2}\mathcal{L}_{II}S=\sum_{i,j,k}h_{ijk}^2-S^2,
\\\frac{1}{2}\mathcal{L}_{II}H^2=|\nabla H|^2-H^2S,
\end{eqnarray}
where $h_{ijk}$ is the component of the covariant derivative of the second fundamental form.
Therefore by the Gauss equation (\ref{gauss}) we have
\begin{eqnarray}\label{T1}
\frac{1}{2}\mathcal{L}_{II}R=|\nabla H|^2-\sum_{i,j,k}h_{ijk}^2-RS.
\end{eqnarray}
\subsection{Self-expanders}
For self-expanders in a Euclidean space, one can introduce the following elliptic operator
\begin{eqnarray}
\mathcal{L}_{III}f=\Delta f+\langle X, \nabla f\rangle,
\end{eqnarray}
where $\Delta$ and $\nabla$ denote the Laplacian and the gradient operator on the self-expander, respectively
and $\langle\cdot,\cdot\rangle$ denotes the standard inner product of $\mathbb{R}^{n+1}$.

We have(cf. \cite{CZ})
\begin{eqnarray}
\frac{1}{2}\mathcal{L}_{III}S&=&\sum_{i,j,k}h_{ijk}^2-(S+1)S,\\
\frac{1}{2}\mathcal{L}_{III}H^2&=&|\nabla H|^2-(S+1)H^2,
\end{eqnarray}
where $h_{ijk}$ is the component of the covariant derivative of the second fundamental form. Then by the Gauss equation (\ref{gauss}) we have
\begin{eqnarray}\label{SE1}
\frac{1}{2}\mathcal{L}_{III}R=|\nabla H|^2-\sum_{i,j,k}h_{ijk}^2-R(S+1).
\end{eqnarray}

\section{Proofs}

 \textbf{Proof of Theorem \ref{main thm}.} Now assume that $M$ is an $n$-dimensional complete self-shrinker in $\mathbb{R}^{n+1}$ satisfying assumptions of Theorem \ref{main thm}. Then from equation (\ref{gauss}), we see that
 \begin{eqnarray}
 H^2-S=R>0,
 \end{eqnarray}
 and
 \begin{eqnarray}
 H^2|\nabla H|^2=\sum_k(\sum_{ij}h_{ij}h_{ijk})^2\leq\sum_{ij}h^2_{ij}\sum_{i,j,k}h^2_{ijk}=S\sum_{i,j,k}h^2_{ijk},
 \end{eqnarray}
which implies that $$|\nabla H|^2\leq\sum_{i,j,k}h^2_{ijk}.$$
Then from equation (\ref{Main equ3}) we obtain
 $$R(1-S)=\sum_{i,j,k}h^2_{ijk}-|\nabla H|^2\geq0,$$ which implies that
 \begin{eqnarray}\label{S}
 S\leq1.
 \end{eqnarray}
 Therefore from Theorem 1.1 in \cite{CL} we see that $M$ is isometric to $\mathbb{S}^m(\sqrt{m})\times \mathbb{R}^{n-m}, 2\leq m\leq n-1$ or  $\mathbb{S}^n(\sqrt{n})$. \endproof
\\
\\
\textbf{Proof of Theorem \ref{main thm2}.} From equation (\ref{Main equ3}) since $R=0$ we have
\begin{eqnarray}\label{equa1}
|\nabla H|^2=\sum_{i,j,k}h^2_{ijk}.
\end{eqnarray}
On the other hand, from the Gauss equation (\ref{gauss}) we see that
\begin{eqnarray}\label{equa2}
H^2=\sum_{i,j}h^2_{ij}
\end{eqnarray}
and hence we have
\begin{eqnarray}
H^2|\nabla H|^2=\sum_k(\sum_{i,j}h_{ij}h_{ijk})^2.
\end{eqnarray}
Note that since
$$\sum_k(\sum_{i,j}h_{ij}h_{ijk})^2\leq\sum_{ij}h^2_{ij}\sum_{i,j,k}h^2_{ijk},$$
 with equality hold at non-totally geodesic points if and only if $$h_{ijk}=\lambda_kh_{ij},$$ we get
  $$H^2|\nabla H|^2\leq\sum_{ij}h^2_{ij}\sum_{i,j,k}h^2_{ijk},$$
   with equality hold at non-totally geodesic points if and only if $$h_{ijk}=\lambda_kh_{ij}.$$ Then from (\ref{equa1}) and (\ref{equa2}) we obtain that
   \begin{eqnarray}
   h_{ijk}=\lambda_kh_{ij}
   \end{eqnarray}
  at non-totally geodesic points, where $\lambda_k$ are functions on $M$. Then from the Codazzi equations (\ref{Codazzi}) we see that
  $$\lambda_kh_{ij}=\lambda_ih_{jk}$$
  at non-totally geodesic points.

At a non-totally geodesic point $p\in M$ assume that $h_{ij}=\mu_i\delta_{ij}$, where $\mu_1,...,\mu_n$ are principle curvatures of $M$ at $p$, we have at $p$
$$\lambda_i\mu_j=0, \s if \s i\neq j.$$
If there exists some $\lambda_i\neq 0$, say $\lambda_1\neq0$, then $\mu_j=0, j=2,...,n$ and therefore the sectional curvature of $M$ at $p$ is zero. Therefore $M$ is a complete hypersurface in $\mathbb{R}^{n+1}$ with zero sectional curvature everywhere. Then by Hartman and Nirenberg (cf. \cite{HN}, Theorem III on page 912), we see that $M$ is isometric to $\Gamma\times\mathbb{R}^{n-1}$, where $\Gamma$ is a complete curve in $\mathbb{R}^2$. Then $\Gamma$ is a complete self-shrinker in $\mathbb{R}^2$. \endproof
\\
\textbf{Proof of Theorem \ref{main thm2'}.} By Theorem \ref{main thm2} we see that $M$ is isometric to $\Gamma\times\mathbb{R}^{n-1}$, where $\Gamma$ is a complete self-shrinker in $\mathbb{R}^2$. Assume that $X$ is the position vector of $\Gamma$, then from equation (\ref{eqmain}) we see that $$H=ce^{\frac{|X|^2}{2}},$$ where $c$ is a constant which could be assumed to be nonnegative. If $c=0$, $\Gamma$ is isometric to $\mathbb{R}$ and $M$ is isometric to $\mathbb{R}^n$. If $c>0$, $\Gamma$ is a bounded curve in $\mathbb{R}^2$ with polynomial volume growth, which must be closed. Then $\Gamma$ is isometric to $\mathbb{S}^1$ or an Abresch-Langer curve. This completes the proof of Theorem \ref{main thm2'}. \endproof

\textbf{Proof of Theorem \ref{main thm2.1}.} Note that this theorem was first proved in \cite{MSS}. Here we give a slight different proof of it.

From equation (\ref{T1}), since $M$ has nonnegative constant scalar curvature, we have
$$|\nabla H|^2-\sum_{i,j,k}h^2_{ijk}=RS.$$
First we will show that $R=0$. If $R$ is a positive constant, then from the Gauss equation (\ref{gauss}) we have
\begin{eqnarray*}
|H|^2-\sum_{ij}h^2_{ij}&=&R>0,
\\|H|^2|\nabla H|^2&=&\sum_k(\sum_{i,j}h_{ij}h_{ijk})^2\leq\sum_{ij}h^2_{ij}\sum_{i,j,k}h^2_{ijk},
\end{eqnarray*}
which implies that
$$|\nabla H|^2<\sum_{i,j,k}h^2_{ijk}.$$
Therefore we have $$RS<0,$$ which is impossible.

Now we have that $R=0$. Similarly with the proof Theorem \ref{main thm2} we can obtain that $M$ is a complete hypersurface in $\mathbb{R}^{n+1}$ with zero sectional curvature, then by  Hartman and Nirenberg (cf. \cite{HN}, Theorem III on page 912), we see that $M$ is isometric to $\Gamma\times\mathbb{R}^{n-1}$, where $\Gamma$ is a complete curve in $\mathbb{R}^2$. Then $\Gamma$ is a straight line in $\mathbb{R}^2$ or a complete translator in $\mathbb{R}^{2}$.

To be precise, assume that $M$ is a translator satisfying
$$H=V^N,$$
where $V$ is a constant vector in $\mathbb{R}^{n+1}$. Denote by $W$ the orthogonal projection of $V$ into the plane $P$ where $\Gamma$ lies in,  then $\Gamma$ is a straight line in $P$ if $W$ is a zero vector and $\Gamma$ is a translator in $P$ satisfying
$$H^{\Gamma}=W^N,$$
where $H^\Gamma$ denotes the mean curvature vector of $\Gamma$ in $P$ and $W^N$ denotes the normal projection of $W$ onto the normal bundle of $\Gamma$ in $P$, if $W$ is not a zero vector. In the last case $\Gamma$ is a straight line or is the Grim Reaper curve in $\mathbb{R}^2$(see for example \cite{Hal}). \endproof

\textbf{Proof of Theorem \ref{main thm3.1}.} Note that this theorem was first proved in \cite{AC}, here we give a slight different proof of it which is similar with the proof of Theorem \ref{main thm2.1}. But we only sketch some main steps.

First since $M$ has constant scalar curvature from (\ref{SE1}) we have $$|\nabla H|^2-\sum_{i,j,k}h^2_{ijk}=R(S+1).$$
Then similarly with the proof of Theorem \ref{main thm2.1} we see that $R$ can not be a positive constant, i.e. $R=0$. Then we can prove that $M$ has zero sectional curvature and hence is isometric to $\Gamma\times \mathbb{R}^{n-1}$, by Hartman and Nirenberg's classification theorem. Then it is easy to see that $\Gamma$ is a complete self-expander in $\mathbb{R}^2$. Furthermore if $\Gamma$ is a straight line, $M$ is a hyperplane and if not, $M$ is a self-expanding hyperplane. \endproof
\section{Appendix}
In this appendix we give a result to partially confirm the conjecture proposed in the introduction. We have
\begin{pro}\label{thmconj}
Let $X:M\to \mathbb{R}^{n+1}$ be an $n$-dimensional complete self-shrinker in $\mathbb{R}^{n+1}$. If the scalar curvature $R$ of $M$ is a constant larger than $n-2$, it is isometric to the round sphere $\mathbb{S}^n(\sqrt{n})$.
\end{pro}
\proof  Note that we have proved $S\leq1$(cf. (\ref{S})). Since
$$H^2=R+S\geq R>n-2\geq0,$$
we may assume that $H>0$ and by Newtons's inequality $$H\geq\sqrt{\frac{nR}{n-1}}.$$ Assume that $\lambda_1\geq\cdot\cdot\cdot\geq\lambda_n$ are principle curvatures of $M$, then
\begin{eqnarray*}
H^2-R&=&S
\\&=&\lambda_1^2+\cdot\cdot\cdot+\lambda_n^2\\
&\geq&\frac{1}{n-1}(H-\lambda_n)^2+\lambda^2_n
\\&=&\frac{1}{n-1}H^2-\frac{2H}{n-1}\lambda_n+\frac{n}{n-1}\lambda^2_n.
\end{eqnarray*}
Therefore
\begin{eqnarray*}
\lambda_n&\geq&\frac{H}{n}-\sqrt{(\frac{n-1}{n}H)^2-\frac{n-1}{n}R}
\\&=&\frac{\frac{n-1}{n}R-\frac{n-2}{n}H^2}{\frac{H}{n}+\sqrt{(\frac{n-1}{n}H)^2-\frac{n-1}{n}R}}.
\end{eqnarray*}
Since
\begin{eqnarray*}
&&\frac{n-1}{n}R-\frac{n-2}{n}H^2
\\ &\geq&\frac{n-1}{n}R-\frac{n-2}{n}(1+R)
\\&=&\frac{R-(n-2)}{n},
\end{eqnarray*}
we see that when $R>n-2$, there is a positive constant $\epsilon_n$ which depends only on $n$ such that
$$\lambda_n\geq\frac{R-(n-2)}{H+\sqrt{(n-1)^2H^2-n(n-1)R}}\geq\epsilon_n.$$
Then by Bonnet-Myers' theorem, $M$ is compact and hence the round sphere $\mathbb{S}^n(\sqrt{n})$, by Corollary 3.2 in \cite{Guo}.
\endproof
\\
\\
\textbf{Acknowledgments} This project is inspired by an online report given by professor Qingming Cheng. The first author would like to thank professor Qingming Cheng and professor Zhen Guo for answering his questions related to this project. The first author was partially supported by the NSF of China(No.11501421, No.11771339).

{}
\vspace{1cm}\sc

Yong Luo

Mathematical Science Research Center of Mathematics,

Chongqing University of Technology,

Chongqing, 400054, China

{\tt yongluo-math@cqut.edu.cn}

\vspace{1cm}

Linlin Sun

School of Mathematics and Statistics and Hubei Key Laboratory of Computational Science,

Wuhan University,

Wuhan, 430072, China

{\tt sunll@whu.edu.cn}

\vspace{1cm}
Jiabin Yin

School of Mathematical Sciences,

Xiamen University,

Xiamen, 361005, China

{\tt jiabinyin@126.com}


\begin{thebibliography}{2}
\bibitem[AL]{AL} U. Abresch and J. Langer, The normalized curve shortening flow and homothetic solutions, {J. Diff. Geom.} {\bf23}(1986), 175--196.
\bibitem[AC]{AC} S. Ancari and X. Cheng, Volume properties and rigidity on self-expanders of mean curvature flow, preprint, arXiv:2012.03364.
\bibitem[AV1]{AV1} S. B. Angenent and J. J. L. Velazquez, Asymptotic shape of cusp singularities in curve shortening, {\em Duke Math. J.} {\bf71}(1995), no.1, 71--110.
 \bibitem[AV2]{AV2} S. B. Angenent and J. J. L. Velazquez, Degenerate neckpinches in mean curvature flow, {\em Crelles J. Math.} {\bf482}(1997), 15--66.
 \bibitem[BM]{BM} T. Begley and K. Moore, On short time existence of Lagrangian
mean curvature flow, {\em Math. Ann.} {\bf367}(2017), no.3-4, 1473--1515.
\bibitem[CL]{CL} H. D. Cao and H. Z. Li, A gap theorem for self-shrinkers of the mean curvature flow in arbitrary codimension, {\em Calc. Var. Partial Differential Equations} {\bf46}(2013), no.3-4,879--889.
\bibitem[Che]{Che} Q. M. Cheng, Complete hypersurfaces in a Euclidean space $\mathbb{R}^{n+1}$ with constant scalar curvature, {\em  Indiana Univ. Math. J.} {\bf51}(2002), no.1, 53--68.
    \bibitem[CLW]{CLW} Q. M. Cheng, Z. Li and G. X. Wei, Complete self-shrinkers with constant norm of the second fundamental form, 	 arXiv:2003.11464, preprint.
\bibitem[CO]{CO} Q. M. Cheng and S. Ogata, 2-Dimensional complete self-shrinkers in $\mathbb{R}^3$, {\em Math. Z.} {\bf284}(2016), 537--542.
\bibitem[CW]{CW} Q. M. Cheng and G. X. Wei, A gap theorem of self-shrikers, {\em Trans. Amer. Math. Soc.} {\bf367}(2015), no.7, 4895--4915.
\bibitem[CY]{CY} S. Y. Cheng and S. T. Yau, Hypersurfaces with constant scalar curvature, {\em Math. Ann.} {\bf225}(1977), no.3, 195--204.
\bibitem[CZ]{CZ} X. Cheng and D. T. Zhou, Spectral properties and rigidity for self-expanding solutions of the mean curvature flows, {\em Math. Ann.} {\bf371}(2018), no.1-2, 371--389.
\bibitem[CM]{CM} T. H. Colding and W. P. Minicozzi, Generic mean curvature flow I; Generic singularities, {\em Ann. of Math.} {\bf 175}(2012), 755--833.
\bibitem[DX]{DX} Q. Ding and Y. L. Xin, The rigidity theorems of self-shrinkers, {\em Trans. Amer. Math. Soc.} {\bf366}(2014), no.10, 5067--5085.
\bibitem[EH]{EH} K. Ecker and G. Huisken, Mean curvature evolution of entire
graphs, {\em Ann. of Math. (2)} {\bf130}(1989), no.3, 453--471.
\bibitem[GT]{GT} J. Q. Ge and Z. Z. Tang, Chern conjecture and isoparametric hypersurfaces, Differential geometry, 49--60, {\em Adv. Lect. Math. (ALM),} {\bf22}, Int. Press, Somerville, MA, 2012.
\bibitem[Guo]{Guo} Z. Guo, Scalar curvature of self-shrinkers, {\em J. Math. Soc. Japan} {\bf70}(2018), no.3, 1103--1110.
\bibitem[Hal]{Hal} H. P. Halldorsson, Self-similar solutions to the curve shortening flow, {\em Trans. Amer. Math. Soc.}{\bf 364}(2012), no.10, 5285--5309.
\bibitem[HN]{HN} P. Hartman and L. Nirenberg, On spherical image maps whose Jacobians do not change sign, {\em Amer. J. Math.} {\bf81}(1959), 901--920.
\bibitem[Hui1]{Hui1} G. Huisken, Asympototic behavior for singularities of the mean curvature flow, {\em J. Diff. Geom.} {\bf31}(1990), no.1, 285--299.
\bibitem[Hui2]{Hui2} G. Huisken, Local and global behavior of hypersurfaces moving by mean curvature, {\em Differential geometry, partial differential equations on manifolds(Los Angeles, CA, 1990), Proc. Sympos. Pure Math.54. Part1.,} Amer. Math. Soc., Providence, RI, 1993, 175--191.
\bibitem[HS]{HS} G. Huisken and C. Sinestrari, Convescity estimates formean curvature flowand singularities of mean convex
surfces, {\em Acta Math.} {\bf183}(1999), 45--70.
\bibitem[Ish]{Ish} N. Ishimura, Curvature evolution of plane curves with prescribed opening angle, {\em Bull. Austral.
Math. Soc.} {\bf52}(1995), no.2, 287--296.
\bibitem[LS]{LS} N. Q. Le and N. Sesum, Blow-up rate of the mean curvature during the mean curvature flow and a gap theorem for self-shrinkers, {\em Comm. Anal. Geom.} {\bf19}(2011), no.4, 633--659.
    \bibitem[LXX]{LXX} L. Lei, H. W. Xu and Z. Y. Xu, A new pinching theorem for complete self-shrinkers and its generalization, {\em Sci. China Math.} {\bf63}(2020), no.6, 1139--1152.
    \bibitem[MSS]{MSS} F. Mart$\acute{i}$n, A. Savas-Halilaj and K. Smoczyk, On the topology of translating solitons of the mean
curvature flow, {\em Calc. Var. Partial Differential Equations} {\bf54}(2015), no.3, 2853--2882.
\bibitem[Ros]{Ros} A. Ros, Compact hypersurfaces with constant scalar curvature and a congruence theorem, {J. Diff. Geom.} {\bf27}(1988), 215--220.
        \bibitem[SWY]{SWY} M. Scherfner, S. Weiss and S. T. Yau, A review of the Chern conjecture for isoparametric hypersurfaces in spheres, Advances in geometric analysis, 175--187, {\em Adv. Lect. Math. (ALM)}, {\bf21}, Int. Press, Somerville, MA, 2012.
\bibitem[Sm1]{Sm1} K. Smoczyk, Self-shrinkers of the mean curvature flow in arbitrary
codimension, {\em Int. Math. Res. Not.} {\bf48}(2005), 2983--3004.
\bibitem[Sm2]{Sm2} K. Smoczyk, Self-expanders of the mean curvature flow, {\em  Vietnam J. Math.}, online, https://doi.org/s10013-020-00469-1.
\bibitem[WXZ]{WXZ} H. J. Wang, H.W. Xu and E. T. Zhao, A global pinching theorem for complete translating
solitons of mean curvature flow, {\em Pure Appl. Math. Q. } {\bf12}(2016), no.4, 603--619.
\bibitem[Wh1]{Wh1} B. White, The size of the singular sets in mean curvature flow of mean convex sets, {\em J. Amer. Math. Soc.} {\bf13}(2000), no.3, 665--695.
\bibitem[Wh2]{Wh2} B. White, The nature of singularities in mean curvature flow of mean convex sets, {\em J. Amer. Math. Soc.} {\bf16}(2003), no.1, 123--138.
\bibitem[Xin]{Xin} Y. L. Xin, Translating solitons of the mean curvature flow, {\em  Calc. Var. Partial Differential Equations} {\bf54}(2015), no.2, 1995--2016.
\bibitem[XX]{XX} H. W. Xu and Z. Y. Xu,  On Chern's conjecture for minimal hypersurfaces and rigidity of self-shrinkers, {\em  J. Funct. Anal.} {\bf273}(2017), no.11, 3406--3425.
\end{thebibliography}
\end{document}